\documentclass[12pt, a4paper]{amsart}
\usepackage{tikz-cd}
\usepackage{amssymb,amsmath}
\usepackage{array}

\usepackage{amsmath, amssymb, amsfonts, amstext, amsthm, amscd}
\usepackage[mathscr]{euscript}
\usepackage{enumerate}

\usepackage{float}
\usepackage{graphicx}

\usepackage[raiselinks=false,colorlinks=true,citecolor=blue,urlcolor=blue,linkcolor=blue,bookmarksopen=true,pdftex]{hyperref}

\usepackage{tikz}
\usepackage{tikz-cd}
\usetikzlibrary{snakes}
\usetikzlibrary{intersections, calc}

\usepackage[english]{babel}
\usepackage{chngcntr}

\begin{document}

\thispagestyle{empty}

\def\theequation{\arabic{section}.\arabic{equation}}

\newcommand{\codim}{\mbox{{\rm codim}$\,$}}
\newcommand{\stab}{\mbox{{\rm stab}$\,$}}
\newcommand{\lr}{\mbox{$\longrightarrow$}}

\newcommand{\be}{\begin{equation}}
\newcommand{\ee}{\end{equation}}

\newtheorem{guess}{Theorem}[section]
\newcommand{\bth}{\begin{guess}$\!\!\!${\bf }~}
\newcommand{\eeth}{\end{guess}}
\renewcommand{\bar}{\overline}
\newtheorem{propo}[guess]{Proposition}
\newcommand{\bpropo}{\begin{propo}$\!\!\!${\bf }~}
\newcommand{\epropo}{\end{propo}}

\newtheorem{lema}[guess]{Lemma}
\newcommand{\blem}{\begin{lema}$\!\!\!${\bf }~}
\newcommand{\elem}{\end{lema}}

\newtheorem{defe}[guess]{Definition}
\newcommand{\bdefe}{\begin{defe}$\!\!\!${\bf }~}
\newcommand{\edefe}{\end{defe}}

\newtheorem{coro}[guess]{Corollary}
\newcommand{\bcor}{\begin{coro}$\!\!\!${\bf }~}
\newcommand{\ecor}{\end{coro}}

\newtheorem{rema}[guess]{Remark}
\newcommand{\brem}{\begin{rema}$\!\!\!${\bf }~\rm}
\newcommand{\erem}{\end{rema}}

\newtheorem{exam}[guess]{Example}
\newcommand{\beg}{\begin{exam}$\!\!\!${\bf }~\rm}
\newcommand{\eeg}{\end{exam}}

\newtheorem{notn}[guess]{Notation}
\newcommand{\bnot}{\begin{notn}$\!\!\!${\bf }~\rm}
\newcommand{\enot}{\end{notn}}
\newcommand{\cw} {{\mathcal W}}
\newcommand{\ct} {{\mathcal T}}
\newcommand{\ch}{{\mathcal H}}
\newcommand{\cf}{{\mathcal F}}
\newcommand{\cd}{{\mathcal D}}
\newcommand{\cR}{{\mathcal R}}
\newcommand{\cv}{{\mathcal V}}
\newcommand{\cn}{{\mathcal N}}
\newcommand{\lra}{\longrightarrow}
\newcommand{\ra}{\rightarrow}
\newcommand{\blr}{\Big \longrightarrow}
\newcommand{\da}{\Big \downarrow}
\newcommand{\ua}{\Big \uparrow}
\newcommand{\hra}{\mbox{{$\hookrightarrow$}}}
\newcommand{\rt}{\mbox{\Large{$\rightarrowtail$}}}
\newcommand{\dua}{\begin{array}[t]{c}
\Big\uparrow \\ [-4mm]
\scriptscriptstyle \wedge \end{array}}
\newcommand{\ctext}[1]{\makebox(0,0){#1}}
\setlength{\unitlength}{0.1mm}
\newcommand{\cm}{{\mathcal M}}
\newcommand{\cl}{{\mathcal L}}
\newcommand{\cp}{{\mathcal P}}
\newcommand{\ci}{{\mathcal I}}
\newcommand{\bz}{\mathbb{Z}}
\newcommand{\cs}{{\mathcal s}}
\newcommand{\ce}{{\mathcal E}}
\newcommand{\ck}{{\mathcal K}}
\newcommand{\cz}{{\mathcal Z}}
\newcommand{\cg}{{\mathcal G}}
\newcommand{\cj}{{\mathcal J}}
\newcommand{\cc}{{\mathcal C}}
\newcommand{\ca}{{\mathcal A}}
\newcommand{\cb}{{\mathcal B}}
\newcommand{\cx}{{\mathcal X}}
\newcommand{\co}{{\mathcal O}}
\newcommand{\bq}{\mathbb{Q}}
\newcommand{\bt}{\mathbb{T}}
\newcommand{\bh}{\mathbb{H}}
\newcommand{\br}{\mathbb{R}}
\newcommand{\bl}{\mathbf{L}}
\newcommand{\wt}{\widetilde}
\newcommand{\im}{{\rm Im}\,}
\newcommand{\bc}{\mathbb{C}}
\newcommand{\bp}{\mathbb{P}}
\newcommand{\ba}{\mathbb{A}}
\newcommand{\spin}{{\rm Spin}\,}
\newcommand{\ds}{\displaystyle}
\newcommand{\tor}{{\rm Tor}\,}
\newcommand{\bff}{{\bf F}}
\newcommand{\bs}{\mathbb{S}}
\def\ns{\mathop{\lr}}
\def\nssup{\mathop{\lr\,sup}}
\def\nsinf{\mathop{\lr\,inf}}
\renewcommand{\phi}{\varphi}
\newcommand{\tT}{{\widetilde{T}}}
\newcommand{\tG}{{\widetilde{G}}}
\newcommand{\tB}{{\widetilde{B}}}
\newcommand{\tC}{{\widetilde{C}}}
\newcommand{\tW}{{\widetilde{W}}}
\newcommand{\tphi}{{\widetilde{\Phi}}}
\newcommand{\fp}{{\mathfrak{p}}}

\title[Equivariant $K$-theory of flag Bott manifolds]{Equivariant $K$-theory of flag Bott manifolds of general Lie
  type}

\author[B. Paul]{Bidhan Paul}
\address{Department of Mathematics, Indian Institute of Technology, Madras, Chennai 600036, India}
\email{bidhan@smail.iitm.ac.in}

\author[V. Uma]{Vikraman Uma}
\address{Department of Mathematics, Indian Institute of Technology, Madras, Chennai 600036, India
}
\email{vuma@iitm.ac.in}



\begin{abstract}
  The aim of this paper is to describe the equivariant and ordinary
  Grothendieck ring and the equivariant and ordinary topological
  $K$-ring of flag Bott manifolds of general Lie type. This will
  generalize the results on the equivariant and ordinary cohomology of
  flag Bott manifolds of general Lie type due to Kaji Kuroki Lee and
  Suh.
\end{abstract}

\maketitle

\section{Introduction}\label{flag}
A Bott tower $M_{\bullet}=\{M_j\mid 1\leq j\leq r\}$, where $M_{j}=\mathbb{P}(1\oplus L_j)$ for a line bundle $L_j$ over $M_{j-1}$ for every $1\leq j$, is a sequence of $\mathbb{P}^1_{\mathbb{C}}$-bundles. Here $M_0=\star$ is a point. At each stage $M_j$ is a $j$-dimensional nonsingular toric variety with an action of a dense torus $(\mathbb{C}^*)^{j}$. We call $M_j$ a $j$-stage Bott manifold for $1\leq j\leq r$.

The topology and geometry of these varieties have been widely studied
recently (see \cite{cms}). One main motivation to study these
varieties comes from its relation to the Bott Samelson variety which
arise as desingularizations of Schubert varieties in the full flag
variety(see \cite{bs}, \cite{d}, \cite{h}). Indeed it has been shown
by Grossberg and Karshon (see \cite{gk}) that a Bott Samelson variety
admits a degeneration of complex structures with special fibre a Bott
manifold (also see \cite{p}). Thus the topology and geometry of a Bott
manifold gets related to those of the Bott -Samelson variety which in
turn is equipped with rich connections with representation theory of
semisimple Lie groups.

The construction of a Bott tower can be generalized in several
directions. One such natural direction is to replace
$\mathbb{P}^1_{\mathbb{C}}$ with complex projective spaces of
arbitrary dimensions. In other words at each stage $M_j$ is the
projectivization of direct sum of finitely many complex line
bundles. In this way we can construct the so called generalized Bott
manifolds and the genralized Bott tower which also have the structure
of a nonsingular toric variety and have been studied widely (see
\cite{cms1}, \cite{kl}).

The construction of a Bott tower has recently been generalized in
another natural direction of a {\it flag Bott manifold} by Kaji,
Kuroki, Lee, Song and Suh (see \cite{kkls}, \cite{klss}, \cite{kur})
by replacing $\mathbb{P}^1_{\mathbb{C}}$ (which can be identified with
the variety of full flags in $\mathbb{C}^2$) at every stage with the
variety $Flag(\mathbb{C}^n)$ of full flags in $\mathbb{C}^n$ for any
$n\geq 2$. More precisely, at the $j$th stage we let $M_j$ to be the
flagification of direct sum of $n_j+1$ complex line bundles over
$M_{j-1}$. Thus $M_j$ is a flag bundle over $M_{j-1}$ with fibre the
full flag manifold $Flag(\mathbb{C}^{n_j+1})$ of dimension
$(n_j)(n_j+1)/2$.

Moreover, in \cite{fls} the {\it flag Bott Samelson variety} has been
introduced and studied by Fujita, Lee and Suh. The flag Bott Samelson
varieties are special cases of the more general notion of the
generalized Bott-Samelson variety which was introduced by Jantzen
\cite{j} and studied by Perrin \cite{per} to obtain small resolutions
of Schubert varieties. Earlier in \cite{sv}, Bott-Samelson varieties
have been used by Sankaran and Vanchinathan to obtain resolutions of
Schubert varieties in a partial flag variety $G/P$ by generalizing
Zelevinsky's results who constructed small resolutions for all
Schubert varieties in Grassmann varieties \cite{z}. The special
property of the flag Bott Samelson variety studied in \cite{fls} is
that they are iterated bundles with fibres full flag varieties. In
particular, similar to the Bott-Samelson variety the flag
Bott-Samelson variety admits a one parameter family of complex
structures which degenerates to a flag Bott manifold(see \cite[Section
4]{fls}). Thus the geometry and topology of flag Bott Samelson variety
gets related with that of flag Bott manifolds since the underlying
differentiable structure is preserved under the deformation and can
therefore be studied using this degeneration.

There is a natural effective action of a complex torus $\mathbb{T}$
and the corresponding compact torus
$\mathbb{T}_{comp}\subseteq \mathbb{T}$ on the flag Bott manifold (see
\cite{kkls}, \cite[Section 3.1]{klss} and Section \ref{torusaction} below).

The flag Bott manifolds of general Lie type has been defined by Kaji,
Kuroki, Lee and Suh in \cite{kkls}, where the
$\mathbb{T}_{comp}$-equivariant cohomology of flag Bott manifolds have
been described.

In \cite{cr} Civan and Ray gave presentations for any complex oriented
cohomology ring (in particular the $K$-ring) of a Bott tower. They also
determine the $KO$-ring for several families of Bott towers.

In \cite{su2} P. Sankaran and the second named author describe the
structure of the topological $K$-ring of a Bott tower, the topological
$K$-ring of a Bott Samelson variety as well as the Grothendieck ring
of a Bott Samelson variety.

In \cite{w1,w2} Mathieu Willems has studied the equivariant and
ordinary $K$-ring of Bott towers and Bott Samelson varieties with
application to Schubert calculus.

In \cite{pu} we study the topological $K$-theory of flag Bott
manifolds of type $A$. More precisely, by applying the classical
results in \cite{at} and \cite{k}, for $K$-ring of a flag bundle
iteratively we give a presentation in terms of generators and
relations of the topological $K$-ring and hence the Grothendieck ring
of the flag Bott manifolds and the flag Bott Samelson variety of type
$A$ (see \cite[Theorem 4.4, Corollary 4.7, Corollary 4.8]{pu}). We
also show that our presentation generalizes the presentation of the
$K$-ring of Bott manifolds and the Bott-Samelson variety in
\cite[Theorem 5.3, Theorem 5.4]{su2}.

In this article we give an alternate algebraic geometric construction
(see Definition \ref{genliebott} below) of flag Bott manifolds of
general Lie type motivated by the results in \cite[Section 1]{u}. It
can be seen that the underlying differentiable manifold structure of
our construction agrees with that of the flag Bott manifold defined in
\cite[Definition 3.5]{kkls}.

Our main aim in this paper is to study the
$\mathbb{T}_{comp}$-equivariant as well as the ordinary topological
$K$-ring of flag Bott manifolds of general Lie type. In Theorem
\ref{kringgenfbm} (resp. Theorem \ref{equivkringgenfbm}) we describe
the ordinary (resp. $\mathbb{T}_{comp}$-equivariant) topological
$K$-ring. In the particular case of flag Bott manifolds of type $A$,
in Theorem \ref{prestypeAfbm}, we derive the presentation given by the
authors in \cite[Theorem 4.4]{pu} by alternate methods. We also give
the presentation for the $\mathbb{T}_{comp}$-equivariant topological
$K$-ring of flag Bott manifolds of type $A$ in Theorem
\ref{equivprestypeAfbm}. This generalizes the results of \cite{kkls}
on the equivariant cohomology to topological $K$-theory. This also
generalizes the results by the authors in \cite[Section 4.1]{pu} to
the equivariant setup.

The flag Bott Samelson variety of general Lie type degenerates to a
flag Bott manifold of general Lie type \cite[Proposition 4.6,
Corollary 4.7]{fls}, where the degeneration preserves the underlying
differentiable structure. Thus in particular we get a presentation for
the equivariant as well as the ordinary $K$-ring of a flag Bott
Samelson variety of general Lie type thereby generalizing our results
in \cite[Section 4.2]{pu}.

Since by construction, the flag Bott manifolds are iterated
$\mathbb{T}$-equivariant cellular fibrations, the fibers at each stage
being generalized flag varieties which have an algebraic cellular
decomposition given by the generalized Schubert cells. Thus by the
equivariant cellular fibration lemma (see \cite[Section 5.5,
Proposition 5.5.6]{cg}) it follows that the ordinary
(resp. $\mathbb{T}$-equivariant) Grothendieck ring of flag Bott
manifolds is isomorphic to the ordinary (resp.
$\mathbb{T}_{comp}$-equivariant) topological $K$-ring. Thus we also
get a description of the $\mathbb{T}$-equivariant as well as the
ordinary Grothendieck ring of flag Bott manifolds.

\section{$K$-theory of flag Bott manifolds of general Lie type}\label{genlie}
 
The flag Bott manifold has been defined for other Lie types in
\cite{kkls} and their $\mathbb{T}_{comp}$-equivariant cohomology ring
has been computed.  In this section we generalize their results to the
topologial $K$-ring as well as the Grothendieck ring of flag Bott
manifold of general Lie type. We give below an alternate algebraic
geometric construction of flag Bott manifolds by referring to the
results in \cite{u}.

Let $G$ be a complex semisimple simply connected algebraic group and
$B$ a Borel subgroup of $G$, $T\subseteq B$ a maximal torus and
$W=N(T)/T$ be the Weyl group.  Let
$\mbox{rank}(G)=\mbox{dim}(T)=n$. Let $P\supseteq B$ denote a
parabolic subgroup of $G$. Let $R_u(P)$ be the unipotent radical of
$P$ and $L=P/R_u(P)$ denote the Levi subgroup of $P$. Let $W_{L}$
denote the parabolic subgroup of $W$ associated to $P$.

Let $K$ (resp.  $T_{comp}=K\cap T$) denote the maximal compact
subgroup of $G$ (resp. $T$). Then $B\cap K=T\cap K=T_{comp}$ is the
maximal torus of $G$ and $P\cap K=L\cap K=Z$ is the centralizer of a
circle subgroup of $T_{comp}$.

The homogeneous space $G/P$ is a flag manifold of general Lie type. In
the case when $G=SL_{n+1}(\mathbb{C})$ and $P=B$ is the subgroup of
upper triangular matrices of determinant $1$, $G/P$ is the complete
flag manifold $Flag(\mathbb{C}^{n+1})$ in type $A$.

Let $H$ be a connected complex linear algebraic group. Let $Y$ be a
nonsingular complex algebraic variety. Let $E$ be a nonsingular
complex algebraic variety with a right action of $H$ such that
\[{E}\lra {Y=E/H}\] is a principal $H$-bundle.

Let $\Phi:H\lra T$ be a morphism of algebraic groups where
$T:=(\mathbb{C}^*)^n$ is an algebraic torus. Note that $\Phi$ induces
$\Phi^*:R(T)\lra R(H)$ by pull back. For, any $T$ representation can be
viewed as a $H$-representation through the morphism $\Phi$.

We recall below the definition of a $T$-cellular variety (see
\cite[Example 1.9.1]{ful}).  \bdefe\label{gcell} A {\em $T$-cellular
  variety} is a $T$-variety $X$ equipped with a $T$-stable algebraic
cell decomposition. In other words there is a filtration
\[X=X_1\supseteq X_2\supseteq\cdots\supseteq X_m\supseteq X_{m+1}=\emptyset\] where
each $X_i$ is a closed $T$-stable subvariety of $X$ and
$X_i\setminus X_{i+1}=U_i$ is $T$-equivariantly isomorphic to the
affine space ${\mathbb{C}}^{k_i}$ equipped with a linear action of
$T$, for $1\leq i\leq m$. Let $x_i\in U_i$ for $1\leq i\leq m$ be the
$T$-fixed points of $X$. \edefe

Let $E(X):=E\times_{H} X$ denote the bundle associated to
$E\lra E/H=Y$ with fiber the $T$-cellular variety $X$ and base
$Y$. Here $H$ acts on $X$ via the map $\Phi:H\lra T$.  Then it is
known (see \cite[Proposition 23]{eg}) that $E(X)$ is a complex
algebraic variety. Further, when $X$ is nonsingular, $E(X)$ is also
nonsingular (see \cite{u}). Let $\pi:E(X)\lra Y$ denote the bundle
projection.

Recall that for a nonsingular $T$-variety $X$ the Grothendieck group
of $T$-equivariant coherent sheaves on $X$ can be identified with the
Grothendieck ring of equivariant vector bundles on $X$ denoted by
$\mathscr{K}_{T}^0(X)$. Similarly the ordinary Grothendieck group of
coherent sheaves on $X$ can also be identified with the Grothendieck
ring of vector bundles on $X$ denoted by $\mathscr{K}^0(X)$ (see
\cite{ful}, \cite{cg}). Recall that $\mathscr{K}^0_{T}(pt)=R(T)$ which
is the Grothendieck ring of complex representations of $T$ (recall
that $\mathscr{K}_{T}^1(pt)=0$). Further, the structure morphism
$X\lra \mbox{Spec} ~\mathbb{C}$ induces canonical $R(T)$-module
structure on $\mathscr{K}_{T}^0(X)$.

We have the pull back morphism
$\pi^*:\mathscr{K}^0(Y)\lra \mathscr{K}^0(E(X))$. Also for a
$T$-equivariant vector bundle $\mathcal{W}$ on $X$ we can consider the
vector bundle $E\times_{H} \mathcal{W}$ on $E(X)$ where $H$
acts on $\mathcal{W}$ via the map $\Phi:H\lra T$. This induces a
morphism $\mathscr{K}^0_{T}(X)\lra \mathscr{K}^0(E(X))$ which sends
\[ [\mathcal{W}]\lra [E\times_{H}\mathcal{W}] .\] Thus we can
consider
\be\label{map1} \phi':\mathscr{K}^0(Y)\otimes \mathscr{K}^0_{T}(X)\lra
  \mathscr{K}^0({E}(X))\ee which takes
$[\mathcal{V}]\otimes [\mathcal{W}]$ to
$[\pi^*(\mathcal{V})]\otimes [E\times_{H} \mathcal{W}]$. Consider a
$T$-representation $\mathbf{W}$ which corresponds to a trivial
$T$-equivariant vector bundle $X\times \mathbf{W}$ on $X$. Note that
$E\times_{H} (X\times \mathbf{W})$ on $E\times_{H} X$ is isomorphic to
$\pi^*(E\times_{H} \mathbf{W})$ where $E\times_{H} \mathbf{W}$ is a
vector bundle on $Y=E/H$. Thus in $\mathscr{K}^0({E}(X))$
$\phi'((E\times_{H} \mathbf{W})\otimes 1)=\phi'(1\otimes (X\times
\mathbf{W}))$ for $[\mathbf{W}]\in R(T)$. Thus $\phi'$ induces a
morphism
\be\label{map2}\phi:\mathscr{K}^0(Y)\otimes_{R(T)}\mathscr{K}^0_{T}(X)\lra
  \mathscr{K}^0({E}(X)).\ee

\bth\label{algcellularkunneth} For a nonsingular $T$-cellular variety
$X$ the map $\phi$ induces the following isomorphism of
$\mathscr{K}^0(Y)$-modules:
\be\label{algkunneth}\mathscr{K}^0({E}(X))\simeq
\mathscr{K}^0(Y)\otimes_{R(T)} \mathscr{K}^0_{T}(X)\ee (see
\cite[Corollary 1.8]{u}). Here $\mathscr{K}^0(Y)$ gets the structure
of an $R(T)$-module by sending $[\mathbf{V}]\in R(T)$ to
$[E\times_{H} \mathbf{V}]\in \mathscr{K}^0(Y)$, where a
$T$-representation $\mathbf{V}$ is thought of as a $H$-representation
via the morphism $\Phi$. \eeth
\begin{proof}
  Let $\alpha:X_{i+1}\lra X_i$ denote the closed embedding.
  Let $\beta: U_i\lra X_i$ denote the immersion of the Zariski open
  set. Since $\alpha$ and $\beta$ are maps of $T$-varieties they induce
  maps in higher $T$-equivariant $K$-theory:
  \[ \alpha_{*}: \mathscr{K}^l_{T}(X_{i+1})\lra
    \mathscr{K}^l_{T}(X_i) \] and
  \[\beta^*:\mathscr{K}^l_{T}(X_i)\lra \mathscr{K}^l_{T}(U_i).\] When
  $l=0$ the maps are defined respectively by
  \[ [\mathcal{G}]\mapsto [\alpha_*(\mathcal{G})]\] where
  $\mathcal{G}$ denotes a $T$-equivariant coherent sheaf on $X_{i+1}$
  and \[[\mathcal{V}]\mapsto [\mathcal{V}\mid_{Z_i}]\] where
  $\mathcal{V}$ denotes a $T$-equivariant vector bundle on $X_i$.
  Now, $\alpha_*$ and $\beta^*$ induce long exact sequence in
  $T$-equivariant algebraic $K$-theory (see \cite{cg}, \cite{mer}):
\[
\begin{array}{llllll}
\cdots\lra &\mathscr{K}^1_{T}(X_i)&
                                \stackrel{\beta^*}{\lra}&\mathscr{K}^1_{T}(U_i)&
                                                                                 \stackrel{\partial}{\lra}&
                                                                                                            \mathscr{K}^0_{T}(X_{i+1})\\
  \stackrel{\alpha_*}{\lra} & \mathscr{K}^0_{T}(X_i)&
                                                      \stackrel{\beta^*}{\lra}\mathscr{K}^0_{T}(U_i)&\lra 0. 
\end{array}  \]

Moreover, since \[ \begin{array}{ll} &U_i\\ &\downarrow\\ &pt\end{array}\]
is a $T$-equivariant affine bundle by Thom isomorphism it follows that
$\mathscr{K}_{T}^l(pt)\simeq \mathscr{K}_{T}^l(U_i)$ for $l\geq 0$ (see
\cite{cg}, \cite{mer}). Now, $\mbox{ker} \alpha_*=\mbox{Im} \partial$
and $\mathscr{K}_{T}^1(U_i)=\mathscr{K}^1_{T}(pt)=0$. Thus $\mbox{Im}
\partial=0$ which implies that $\mbox{ker} \alpha_*=0$. Thus we have a
short exact sequence of $R(T)$-modules:
\be\label{ses1}0\ra  \mathscr{K}^0_{T}(X_{i+1})
                       \stackrel{\alpha_*}{\ra}
    \mathscr{K}^0_{T}(X_i)
                              \stackrel{\beta^*}{\ra}\mathscr{K}^0_{T}(U_i)\ra
                              0.\ee

 Moreover, \eqref{ses1} splits since   \be\label{ti}\mathscr{K}^0_{T}(U_i)\simeq
 \mathscr{K}^0_{T}(pt)\simeq R(T)\ee is a free $R(T)$-module. Tensoring
 \eqref{ses1} over $R(T)$ with $\mathscr{K}^0(Y)$ we get a short exact
 sequence of $\mathscr{K}^0(Y)$-modules:

 \be\label{ses2}0\lra 
   \mathscr{K}^0(Y)\otimes_{R(T)}\mathscr{K}^0_{T}(X_{i+1})
   \stackrel{\mbox{Id}\otimes \alpha_*}{\ra} 
   \mathscr{K}^0(Y)\otimes_{R(T)}\mathscr{K}^0_{T}(X_i) \stackrel{\mbox{Id}\otimes
     \beta^*}{\ra}\mathscr{K}^0(Y)\otimes_{R(T)}\mathscr{K}^0_{T}(U_i)\ra 0.\ee

   By \eqref{ti} the last term of the above exact sequence is
   isomorphic to $\mathscr{K}^0(Y)$. Let $E(X_i):=E\times_{H} X_i$
   where $H$ acts on $X_i$ through $\Phi:H\lra T$ for $1\leq i\leq m$.
   Consider
 \[E(X)=E(X_1)\supseteq E(X_2)\supseteq \cdots \supseteq
   E(X_m)=E(x_m)=Y\] a cellular fibration over $Y$. Also
 $E(X_i)\setminus E(X_{i+1})=E(U_i) = E\times_{H} U_i$ where
 \be\label{affbun}\begin{array}{ll}E\times_{H} U_i \\~~~~\da \\ ~~~~Y\end{array}\ee is
 an affine bundle over $Y$. Considering the closed embedding
 $\alpha':E(X_{i+1})\lra E(X_i)$ and Zariski open embedding
 $\beta': E(U_i)\lra E(X_i)$ we get a long exact sequence in algebraic
 $K$-theory:

\[
\begin{array}{llllll}
\cdots\ra \mathscr{K}^1(E(X_i))
                                \stackrel{\beta'^*}{\ra}\mathscr{K}^1(E(U_i))&
                                                                                 \stackrel{\partial}{\ra}&
                                                                                                            \mathscr{K}^0(E(X_{i+1}))\\
  \stackrel{\alpha'_*}{\ra}  \mathscr{K}^0(E(X_i))
                                                      \stackrel{\beta'^*}{\ra}\mathscr{K}^0(E(U_i))&\ra 0. 
\end{array}  \]

We further have the commuting triangle
\begin{center}
\begin{tikzcd}[column sep=scriptsize]
	\mathscr{K}^l(E(X_i)) \arrow[rr, "\beta'^*"blue]
	&
	& \mathscr{K}^l(E(U_i))  \\
	& \mathscr{K}^l(Y)\arrow[ul] \arrow[ur, "\cong"blue]
	&
\end{tikzcd}\end{center}

  where the isomorphism is the Thom isomorphism
       for the affine bundle \eqref{affbun}. Thus it follows that
       $(\beta')^*$ is surjective which in turn implies that the
       connecting homomorphism $\partial$ is trivial and $\alpha'_{*}$
       is injective. Thus we get the following short exact sequence of
       $\mathscr{K}^0(Y)$ modules: \be\label{ses3} 0\lra 
                                                                                                            \mathscr{K}^0(E(X_{i+1}))
  \stackrel{\alpha'_*}{\lra} \mathscr{K}^0(E(X_i))
                                                      \stackrel{\beta'^*}{\lra}\mathscr{K}^0(E(U_i))
                                                      \lra 0.\ee

                                                      From
                                                      \eqref{ses2} and
                                                      \eqref{ses3} we
                                                      get the
                                                      following
                                                      commuting
                                                      diagram:

\begin{center}
	\begin{tikzcd}[cramped,column sep=tiny]
		0\ar[r]	&\mathscr{K}^0(Y)\otimes_{R(T)}\mathscr{K}^0_T(X_{i+1}) \arrow[r, "Id\otimes \alpha_*" red] \arrow[d, "\phi_{i+1}" blue]
		&\mathscr{K}^0(Y)\otimes_{R(T)}\mathscr{K}^0_T(X_{i}) \arrow[r, "Id\otimes \beta^*" red] \arrow[d, "\phi_i" blue] 
		&\mathscr{K}^0(Y)\otimes_{R(T)}\mathscr{K}^0_T(U_{i})  \arrow[d,-,double equal sign distance,double]\ar[r]&0\\
		0\ar[r]	&\mathscr{K}^0(E(X_{i+1})) \arrow[r, "\alpha'_*" red]
		&\mathscr{K}^0(E(X_{i})) \arrow[r, "\beta'^*" red]& \mathscr{K}^0(E(U_{i}))\ar[r]&0\\
	\end{tikzcd}
\end{center}

where the last equality follows from the Thom isomorphism
$\mathscr{K}^0(E(U_i))$ with $\mathscr{K}^0(Y)$. The first downward
arrow $\phi_{i+1}$ is an isomorphism by induction assumption. Thus by
diagram chasing the middle downward arrow $\phi_i$ is an
isomorphism. We can start the induction when $i+1=m$ since $E(x_m)=Y$
and $\mathscr{K}^0_{T}(x_m)=R(T)$. The theorem now follows since
$\phi_1=\phi$ is an isomorphism.
                                                 \end{proof}

Note that $G/P$ has a $T$-stable and hence a $H$-stable ($H$ acts via
$\Phi$) algebraic cell decomposition given by the generalized Schubert
cells $B\cdot w\cdot P/P$ where $w\in W^{L}$ where $W^{L}$ denote the
set of minimal length coset representatives of $W/W_{L}$ (see \cite{kk}).

Consider the associated flag manifold bundle
${E}(G/P):={E}\times_{H} G/P$ where $H$ acts on $G/P$ by composing
the standard left action of $T$ with $\Phi$.

We have the following corollary of Theorem \ref{algcellularkunneth} by
taking $X=G/P$ and $T$ to be the maximal torus of $G$.

\bcor\label{algkunnethflag}
\be\label{algkunnethflagiso}\mathscr{K}^0({E}(G/P))\simeq
\mathscr{K}^0(Y)\otimes_{R(T)} \mathscr{K}^0_{T}(G/P)\ee
\ecor

Let $K^0$ denote the topological $K$-ring (see \cite{at}) and
$K^0_{T_{comp}}$ denote $T_{comp}$-equivariant topological $K$-ring
(see \cite{segal}). Recall that (see \cite{cg}) we have the
isomorphism \be\label{iso3'} R(T)\simeq R(T_{comp}). \ee

Further, for a $T$-cellular variety $X$, $E(X)\lra Y$ is a cellular
fibration in the sense of \cite[Section 5.5, Proposition
5.5.6]{cg}. If the canonical forgetful map
$\mathscr{K}^0(Y)\lra K^0(Y)$ is an isomorphism then we have the
isomorphism: \be\label{iso1}\mathscr{K}^0({E}(X))\simeq K^0({E}(X))\ee
where $E(X):=E\times_{H} X$.

In particular, when $Y=pt$ then
$X\lra pt$ is a $T$-equivariant
cellular fibration where $T$ acts trivially on $pt$. Then by
\cite[Section 5.5, Proposition 5.5.6]{cg} we have: \be\label{iso2}
\mathscr{K}_{T}^0(X)\simeq K_{T_{comp}}^0(X) .\ee

Therefore from \eqref{algkunneth} and \eqref{iso1} and \eqref{iso2} we
have the following theorem. However, we shall give a self-contained
proof of this theorem using the results in topological equivariant
$K$-theory (see \cite{segal}).  In addition we shall assume that for
the space $Y$, $K^{-1}(Y)=0$.

\bth\label{cellularkunnethtop} We have the following isomorphism:
\be\label{topkunneth}{K}^0(E(X))\simeq
{K}^0(Y)\otimes_{R(T_{comp})}{K}^0_{T_{comp}}(X).\ee Here ${K}^0(Y)$
is a $R(T_{comp})=R(T)$-module via the map which takes any
$[\mathbf{V}]\in R(T_{comp})$ for a $T$-representation $\mathbf{V}$ to
$[E\times_{H} \mathbf{V}]\in {K}^0(Y)$, where $H$ acts on $\mathbf{V}$
via $\Phi$.  Further, we have $K^{-1}(E(X))=0$. \eeth
\begin{proof}

By \cite[Proposition 2.6, Definition 2.7, Definition 2.8,
  Proposition 3.5]{segal} it follows that we have a long exact
  sequence of $T_{comp}$-equivariant $K$-groups which is infinite in
  both directions: \begin{align}\label{les} \cdots\ra
  {K}^{-q}_{T_{comp}}(X_i,X_{i+1})\ra {K}_{T_{comp}}^{-q}(X_i)&\ra 
  K^{-q}_{T_{comp}}(X_{i+1})\ra \nonumber\\
                                   &\ra K^{-q+1}_{T_{comp}}(X_i,X_{i+1})\ra 
  \cdots \end{align} for $1\leq i\leq m$ and $q\in \mathbb{Z}$.

Moreover, by \cite[Proposition 2.9]{segal} and \cite[Proposition
3.5]{segal} we have
\begin{align*}K^q_{T_{comp}}(X_i,X_{i+1})
  &=K^q_{T_{comp}}(X_i\setminus X_{i+1})\\ &=
                                             K_{T_{comp}}^q(\mathbb{C}^{k_i})\simeq K^q_{T_{comp}}(x_i)\\
  &=\widetilde{K}^{-q}_{T_{comp}}(x_i^{+})\\
  &=\widetilde{K}_{T_{comp}}^0(S^q)=R(T_{comp})\otimes
    \widetilde{K}^0(S^q(x_i^{+}))\end{align*} for $1\leq i\leq m$ (see
    \cite[Proposition 2.2]{segal}). Thus when $q$ is even
    $K^{-q}_{T_{comp}}(X_i,X_{i+1})=R(T_{comp})$ and when $q$ is odd
    $K^{-q}_{T_{comp}}(X_i,X_{i+1})=0$. Here $x_i^{+}$ is the sum of
    the $T_{comp}$-fixed point $x_i$ and a base point $\mathfrak{o}$
    which is also $T_{comp}$-fixed (see \cite[p. 135]{segal}).


    Alternately, we can also identify
    $K^q_{T_{comp}}(X_i,X_{i+1})={\widetilde
      K^q_{T_{comp}}}(X_i/X_{i+1})$ where
    $X_i/X_{i+1}=Th(U_i)\simeq U_i^{+}\simeq S^{2k_i}$ the Thom space
    of the trivial equivariant vector bundle $U_i$ over the $T$-fixed
    point $x_i$ so that the above isomorphism is the Thom isomorphism
    (see \cite{may}). Indeed for any integer $q$ we have
    $\widetilde{K}^{-q}_{T_{comp}}(S^{2k_i})\simeq
    \widetilde{K}_{T_{comp}}^0(S^{q+2k_i})$ \cite[p. 136]{segal}.

Thus when $q$ is even
$K^{-q}_{T_{comp}}(X_i,X_{i+1})\simeq
\widetilde{K}_{T_{comp}}^0(S^{q+2k_i})=R(T_{comp})$ since $q+2k_i$ is
even. Further, when $q$ is odd
$K^{-q}_{T_{comp}}(X_i,X_{i+1})=\widetilde{K}_{T_{comp}}^0(S^{q+2k_i})=R(T_{comp})\otimes
\widetilde{K}^0(S^{q+2k_i})=0$ since $q+2k_i$ is odd (see
\cite[Section 2, Proposition 3.5, Proposition 2.2]{segal} and
\cite[Lemma 2.7.4]{at}).

Moreover, since $X_m=U_m=\{x_m\}$ and $X_{m+1}=\emptyset$ we have
$K^0_{T_{comp}}(X_m)=R(T_{comp})$ and
$K^{-1}_{T_{comp}}(X_m)=K^{-1}_{T_{comp}}(x_m)=
\widetilde{K}_{T_{comp}}^{-1}(x_m^{+})=\widetilde{K}_{T_{comp}}^0(S^1)=0$
where $x_m^{+}=x_m\sqcup \mathfrak{o}$ where
both $x_m$ and the base point $\mathfrak{o}$ are $T_{comp}$-fixed (see
\cite[p. 135]{segal}).

Now, by decreasing induction on $i$ suppose that
$K_{T_{comp}}^0(X_{i+1})$ is a free $R(T_{comp})$-module for
$1\leq i\leq m$ of rank $m-i$ and $K_{T_{comp}}^{-1}(X_{i+1})=0$.  We
can start the induction since $K^0_{T_{comp}}(X_m)=R(T_{comp})$ and
$K^{-1}_{T_{comp}}(X_m)=0$.

It then follows that from the \eqref{les} we get the following split
short exact sequence of $R(T_{comp})$-modules \be \label{eq1} 0\ra
K_{T_{comp}}^0(X_i, X_{i+1})\ra K_{T_{comp}}^0(X_i)\stackrel
{\alpha^*} {\ra} K_{T_{comp}}^0( X_{i+1})\ra 0 \ee for
$1\leq i\leq m$.

Hence $K_{T_{comp}}(X_i)$ is a free $R(T_{comp})$-module of rank
$m-i+1$. By induction it follows that $K_{T_{comp}}(X_i)$ is free
$R(T_{comp})$-module of rank $m$.

Since $K^{-1}_{T_{comp}}(X_i,X_{i+1})=0$ as shown above and by
induction assumption $K^{-1}_{T_{comp}}(X_{i+1})=0$, it follows from
\eqref{les} that $K^{-1}_{T_{comp}}(X_i)=0$. Thus
$K_{T_{comp}}^{-1}(X_1)=K_{T_{comp}}^{-1}(X)=0$ by induction on $i$.



Now, by tensoring \eqref{eq1} with $K^0(Y)$ over $R(T_{comp})$ we get the
exact sequence

\be\label{eq2} \displaystyle 0\ra K^0(Y)\ra K^0(Y) \otimes_{R(T_{comp})}
K^0_{T_{comp}}(X_i)\ra K^0(Y) \otimes_{R(T_{comp})}
K^0_{T_{comp}}(X_{i+1})\ra 0\ee for $1\leq i\leq m$.

Since $E(X_{i+1})\subseteq E(X_{i})$ is closed subspace for
$1\leq i\leq m$, by \cite[Proposition 2.4.4, Theorem 2.4.9]{at} it
follows that we have a long exact sequence of $K$-groups which is
infinite in both directions: \begin{align}\label{les1} \cdots\ra
  {K}^{-q}(E(X_i),E(X_{i+1}))\ra {K}^{-q}(E(X_i)) &\ra
  K^{-q}(E(X_{i+1}))\ra \nonumber\\ & \ra
  K^{-q+1}(E(X_i),E(X_{i+1}))\ra \cdots \end{align} for
$1\leq i\leq m$.

Thus we have
\begin{align*} K^{-q}(E(X_i),E(X_{i+1}))&=\widetilde{K}^{-q}(E(X_i)/
  E(X_{i+1}))\\ &=\widetilde{K}^{-q}(E(U_i)^+)\\ &\simeq
  K^{-q}(Y)\end{align*} for $1\leq i\leq m$ (see \cite[Definition
2.4.1, Proposition 2.7.2]{at}).



By decreasing induction on $i$ we assume that $K^{-1}(E(X_{i+1}))=0$
and $K^0(E(X_{i+1}))$ is a free $K^0(Y)$-module of rank $m-i$. We can
start the induction since
$K^{-1}(E(X_m))=K^{-1}(E(U_m))=\widetilde{K}^{-1}(E(U_m)^{+})=K^{-1}(Y)=0$
and $K^0(E(X_m))=K^0(E(U_m))=\widetilde{K}^0(E(U_m)^+)=K^0(Y)$. Then
from \eqref{les1} we get the following short exact sequence

\be\label{eq3} 0\lra K^0(Y)\lra K^0(E(X_i))\lra K^0(E(X_{i+1}))\lra
0\ee and we also get that $K^{-1}(E(X_i))=0$ since
$K^{-1}(E(X_i),E(X_{i+1}))=\widetilde{K}^{-1}(E(U_i)^+)=K^{-1}(Y)=0$.

Now, \eqref{eq3} splits as $K^0(Y)$-module so that $K^0(E(X_i))$ is a
free $K^0(Y)$ module of rank $m-i+1$.  By induction we get that
$K^0(E(X))$ is a free $K^0(Y)$-module of rank $m$ and
$K^{-1}(E(X))=0$.


We have the following commuting diagram where the first row is the
exact sequence \eqref{eq2} and the second row is the exact sequence
\eqref{eq3}.
\begin{center}
	\begin{tikzcd}[cramped,column sep=tiny]
		0\ar[r]	
		& K^0(Y) \arrow[r] \arrow[d,-,double equal sign distance,double]
		& K^0(Y)\otimes_{R(T_{comp})}K^0_{T_{comp}}(X_i) \arrow[r] \arrow[d, "\phi_i" blue] 
		&K^0(Y)\otimes_{R(T_{comp})}K^0_{T_{comp}}(X_{i+1})  \arrow[d,  "\phi_{i+1}" blue]\ar[r]
		&0\\
		0\ar[r]	& K^0(Y)\arrow[r]
		&K^0(E(X_i))\arrow[r]
		&K^0(E(X_{i+1}))\ar[r]
		&0\\
	\end{tikzcd}
\end{center}

Also the second vertical arrow is $\phi_i$ and the third
vertical arrow is $\phi_{i+1}$. Since
\[\phi_m:K^0(Y)\otimes_{R(T_{comp})} K^0_{T_{comp}}(x_m=X_m)\lra
K^0(E(x_m))=K^0(Y)\] is the identity map we can start the
induction. Note that the first vertical arrow is identity map and
$\phi_{i+1}$ is an isomorphism by induction hypothesis.

By diagram chasing it follows that $\phi_i$ is an isomorphism

Thus it follows by induction on $i$ that the map
\[\phi_i:K^0(Y) \otimes_{R(T_{comp})} K^0_{T_{comp}}(X_i)\lra
  K^0(E(X_i))\] defined by sending
$[\mathcal{V}]\otimes [\mathcal{W}]$ to
$[\pi^*(\mathcal{V})]\otimes [E\times_{H} \mathcal{W}]$ is an
isomorphism for $1\leq i\leq m$. This proves the theorem since
$\phi_1=\phi$.
\end{proof}  

We have the following topological analogue of Corollary
\ref{algkunnethflag} since $G/P$ is a $T$-cellular variety (see
\cite{kk}).  By putting $X=G/P$ we consider here the $T_{comp}$-action
on $G/P$.

\bcor\label{topkunnethflag}
\be\label{topkunnethflagiso}{K}^0({E}(G/P))\simeq
{K}^0(Y)\otimes_{R(T_{comp})} {K}^0_{T_{comp}}(G/P)\ee
\ecor

We further have the following refinement of \eqref{algkunnethflagiso}
(resp. \eqref{topkunnethflagiso}) using the structure of
$\mathscr{K}^0_{T}(G/P)$ (resp. $K^0_{T_{comp}}(G/P)$).

\bth We have the following isomorphism:

\be\label{algkbun}  \mathscr{K}^0(E(G/P)) \simeq
\mathscr{K}^0(Y)\otimes_{R(T)^{W}} R(T)^{W_L}\ee
\be\label{topkbun} {K}^0(E(G/P)) \simeq {K}^0(Y)\otimes_{R(T_{comp})^{W}}
R(T_{comp})^{W_L}. \ee where $R(T)^{W}$ acts on $\mathscr{K}^0(Y)$
by considering the chain of maps
\[R(T)^{W}\hra R(T)\stackrel{\Phi^*}{ \lra} R(H)\] and taking
$[\mathbf{V}]\in R(T)^{W}$ to the associated bundle
$ [E\times_{H} \mathbf{V}]$ in $K^0(Y)$ by considering $\mathbf{V}$ as
a $H$- representation. \eeth
\begin{proof}
  By \cite[Section 4.1]{mer} (also see \cite{kk}) we have an
  isomorphism \be\label{kk} \mathscr{K}^0_{T}(G/P)\simeq
  R(T)\otimes_{R(T)^W} R(T)^{W_{L}}.\ee
  This follows by identifying
  $R(T)^{W_{L}}$ with $R(P)=\mathscr{K}^0_{G}(G/P)$ and $R(T)^{W}$
  with $R(G)$ and by the isomorphism
  $$R(T)\otimes_{R(G)} \mathscr{K}^0_{G}(G/P)\simeq
  \mathscr{K}^0_{T}(G/P).$$ Thus combining 
  \eqref{algkunnethflagiso} and \eqref{kk} we get
  \be\label{iso3}\mathscr{K}^0(E(G/P))\simeq
  \mathscr{K}^0(Y)\otimes_{R(T)} \otimes (R(T)\otimes_{R(T)^{W}}
  R(T)^{W_{L}}).\ee Since the right hand side of \eqref{iso3} is
  isomorphic to \[\mathscr{K}^0(Y)\otimes_{R(T)^{W}} R(T)^{W_L}\]
  \eqref{algkbun} follows. The proof of \eqref{topkbun}  is
  similar by using \eqref{topkunnethflagiso} and the
  isomorphism
  \[{K}^0_{T_{comp}}(G/P)\simeq R(T_{comp})\otimes_{R(T_{comp})^W}
    R(T_{comp})^{W_{L}}. \]
\end{proof}  

\brem
When $P=B$, $G/P$ is the full flag variety and $W_{L}=W$ then we have
the isomorphism
\[K^0_{T_{comp}}(G/B)\simeq R(T_{comp})\otimes_{R(T_{comp})^W}
  R(T_{comp}).\]  due to McLeod \cite{ml}.
\erem

\subsection{Construction of flag Bott manifolds of general Lie type}

Let $G_j$ be a reductive complex algebraic group of any Lie type,
$B_j\subseteq G_j$ be a Borel subgroup, $T_j\subseteq B_j$ a maximal
torus, $P_j\supseteq B_j$ be a parabolic subgroup and $L_j$ be its Levi
subgroup. Let $W_j:=N(T_j)/T_j$ be the Weyl group for every
$1\leq j\leq r$. Let $\mbox{rank}(G_j)=\mbox{dim}(T_j)=n_j+1$.

Let
$\displaystyle\Phi_k=(\Phi^{(k)}_j)_{j=1}^k: \prod_{j=1}^k P_j\lra
\prod_{j=1}^k L_j$ be the product of the canonical projections
$\Phi^{(k)}_j :P_j\lra L_j $ for $1\leq j\leq k$. Let
$\displaystyle \Psi_{k+1}=(\Psi^{(k+1)}_j)_{j=1}^k\in
Hom_{alg}(\prod_{j=1}^k L_j, T_{k+1})$. Then we have the induced map
$\displaystyle\Psi_{k+1}^*: R(T_{k+1})\ra \bigotimes_{j=1}^k R(L_j)
=R(\prod_{j=1}^k L_j)$. Also, we have
$\displaystyle\Phi_k^*: R(\prod_{j=1}^k L_j)\ra R(\prod_{j=1}^k P_j)$.

We form the following inductive construction:
\be\label{genbm}{M}_{k+1}:= \prod_{j=1}^k G_j\times_{\prod_{j=1}^k
  P_j}G_{k+1}/P_{k+1}\ee where
$\displaystyle \prod_{j=1}^k P_{j}$ acts on $G_{k+1}/P_{k+1}$
via $\Psi_{k+1}\circ \Phi_{k}$ and on
$\displaystyle\prod_{j=1}^k G_j$ by multiplication on the
right as given in \eqref{rightmult}.

Then $M_{k+1}$ is a bundle over
$M_{k}=\prod_{j=1}^k G_j/\prod_{j=1}^k P_j$ with fiber
$G_{k+1}/P_{k+1}$.

\bdefe\label{genliebott}
A flag Bott tower of general Lie type is $\{M_j~:~j=1,\ldots, r\}$
with projections $p_j: M_j\lra M_{j-1}$ for $1\leq j\leq r$ where
$M_j$ is the $j$-stage flag Bott manifold. In particular, $M_0=pt$
and $M_1=G_1/P_1$. Thus we have
\[ M_r\stackrel{p_r}{\lra}\cdots{\lra} M_j\stackrel{p_j}{\lra}
  M_{j-1}\stackrel{p_{j-1}}{\lra} \cdots
  M_1=G_1/P_1\stackrel{p_1}{\lra} M_0=pt. \] \edefe

The above definition is algebraic geometric analogue of the
construction in \cite[Section 3]{kkls} and its underlying
differentiable manifold structure is same as that defined in
\cite[Definition 3.5]{kkls}.

More precisely, let
$(G_{\bullet}, P_{\bullet})=\{(G_j,P_j)\mid 1\leq j\leq r\}$. Given a
family of algebraic group homomorphisms
$\Psi_{j}^{(l)}:L_j\lra T_l\mid 1\leq j<l\leq r\}$, the space $M_r$ is
defined as
$M_r:=G_1\times \cdots\times G_{r}/P_1\times\cdots\times P_{r}$ where
$(p_1,\ldots, p_r)\in P_1\times\cdots\times P_{r}$ acts on
$(g_1,\ldots, g_r)\in G_1\times \cdots\times G_{r}$ from the right by
$(g_1,\ldots, g_r)\cdot (p_1,\ldots, p_r):=$
\be\label{rightmult}(g_1\cdot p_1, \Psi_1^{(2)}\circ \Phi_1^{(r)}
(p_1)^{-1}\cdot g_2\cdot p_2, \prod_{j=1}^2 \Psi_j^{(3)}\circ \Phi_j^{(r)}
(p_j)^{-1}\cdot g_3\cdot p_3, \ldots, \prod_{j=1}^{r-1}
\Psi_j^{(r)}\circ \Phi_j^{(r)} (p_j)^{-1}\cdot g_r\cdot p_r) .\ee This
action is free and hence $M_r$ is a smooth manifold. Moreover, $M_r$
has the structure of a $G_{r}/P_r$ fibre bundle over $M_{r-1}$ (see
\cite{kkls}).

The following proposition describes the Grothendieck ring of $M_r$ as
an algebra over $\mathscr{K}^0(M_{r-1})$ (resp. the topological
$K$-ring of $M_r$ as an algebra over $K^0(M_{r-1})$). See
\cite[Theorem 4.2]{kkls} for the analogous statement for equivariant
cohomology ring.  \bpropo\label{kringgenfbmind} \be\label{grothfbmind}
\mathscr{K}^0(M_r)\simeq \mathscr{K}^0(M_{r-1})\otimes_{R(T_r)^{W_r}}
R(T_r)^{W_{L_{r}}}.\ee \be\label{topfbmind} K^0(M_r)\simeq
K^0(M_{r-1})\otimes_{R(T_r)_{comp}^{W_r}}
R((T_r)_{comp})^{W_{L_{r}}}.\ee Here $\mathscr{K}^0(M_{r-1})$
(resp. $K^0(M_{r-1})$) is an $R(T_r)^{W_r}$
(resp. $R(T_r)_{comp}^{W_r}=R(T_r)^{W_r}$)-algebra via the maps
\[ \displaystyle R(T_r)^{W_r}\hra R(T_r)\stackrel{\Psi_r^*}{\lra}
  \bigotimes_{j=1}^{r-1} R(L_j)=R(\prod_{j=1}^{r-1} L_j)\] where
$\displaystyle[\mathbf{V}]\in R(T_r)^{W_r}$ maps to the class of the
associated vector bundle
$\displaystyle\prod_{j=1}^{r-1}G_j\times_{\prod_{j=1}^{r-1} P_j}
\mathbf{V}$ where $\displaystyle\prod_{j=1}^{r-1} P_j$ acts on
$\mathbf{V}$ via the map $\Psi_r\circ \Phi_{r-1}$. \epropo
\begin{proof}
The isomorphism \eqref{grothfbmind} follows from \eqref{algkbun}
applied when $\displaystyle H=\prod_{j=1}^{r-1} P_j$,
$\displaystyle E=\prod_{j=1}^{r-1} G_j$, $Y=M_{r-1}$ and $T=T_r$. Similarly the
isomorphism \eqref{topfbmind} follows from \eqref{topkbun}.
\end{proof}

Recall (see \cite{st}) that
$R(T_j)=\mathbb{Z}[X^*(T_j)]\simeq\mathbb{Z}[{\bf y}_j]$ where
${\bf y}_j=\{y_{j,1}^{\pm 1},\ldots, y^{\pm 1}_{j,n_j+1}\}$ for
$1\leq j\leq r$.

Then $R(P_j)=R(L_j)=\mathbb{Z}[{\bf y}_j]^{W_{L_j}}=R(T_j)^{{W_{L_j}}}$ for
$1\leq j\leq r$.

 Let
 \[\mathfrak{S}:=\bigotimes_{j=1}^r\mathbb{Z}[{\bf y}_j]^{W_{L_j}}\simeq
   \bigotimes_{j=1}^r R(L_j).\] Let $\mathfrak{I}_j$ denote the ideal
 in $\mathfrak{S}$ generated by the elements
 \[g({\bf y}_j)-\Psi_j^*(g({\bf y}_j))\] where
 $g({\bf y}_j)\in R(T_j)^{W_j}=\mathbb{Z}[{\bf y}_j]^{W_j}\subseteq
 \mathbb{Z}[{\bf y}_j]^{W_{L_j}}$ for every $1\leq j\leq r$. Note that
 $\Psi_j^*(g({\bf y}_j))=g(\Psi_j^* ({\bf y}_j))$ since
 \[\Psi_j^*:R(T_j)=\mathbb{Z}[{\bf y}_j]\lra \bigotimes_{k=1}^{j-1}\mathbb{Z}[{\bf
     y}_{k}]^{W_{L_k}}\subseteq \bigotimes_{j=1}^r\mathbb{Z}[{\bf
     y}_{k}]\] is a $\mathbb{Z}$-algebra homomorphism.)  Also by
 \cite{st} $R(G_j)=R(T_j)^{W_j}$.  Let
 $\mathfrak{I}=\mathfrak{I}_1+\cdots+\mathfrak{I}_r$.

 The following theorem gives a presentation for the toplogical
 $K$-ring of $M_r$. See \cite[Corollary 4.3]{kkls} for the analogous
 statement for equivariant cohomology ring. \bth\label{kringgenfbm} We
 have the following isomorphism of rings
 \be\label{presgenfbm}\mathscr{K}^0(M_r)\simeq K^0(M_r)\simeq
 \mathfrak{S}/\mathfrak{I}.\ee\eeth \begin{proof} The first
   isomorphism follows by \cite[Proposition 5.5.6]{cg} since $M_r$ is
   constructed as iterated cellular fibrations. We now prove the
   second isomorphism in \eqref{presgenfbm}. When $r=1$ then
   $\mathfrak{S}=\mathbb{Z}[{\bf y}_1]^{W_{L_1}}$ and
   $\mathfrak{I}=\mathfrak{I}_1$ is generated by elements
   $g({\bf y}_1)-\Psi_1^*(g({\bf y}_1))$ for
   $g({\bf y}_1)\in R(T_1)^{W_1}$ and $\Psi_1: pt=T_0\lra T_1$. Thus
   $\Psi_1^*: R(T_1)\lra \mathbb{Z}$ is the augmentation map which
   sends $[V]$ to $\mbox{dim}(V)$. But it is known that
   $K^0(G_1/P_1)\simeq
   \mathbb{Z}\otimes_{R(T_1)^{W_1}}R(T_1)^{W_{L_1}}$. Thus the
   statement is true for $r=1$. We assume by induction that the
   statement is true for $r-1$. Thus $K^0(M_{r-1})$ is isomorphic to
   $\mathfrak{S}'=\mathbb{Z}[{\bf y}_1]^{W_{L_1}}\bigotimes \cdots
   \bigotimes\mathbb{Z}[{\bf y}_{r-1}]^{W_{L_{r-1}}}$ modulo the ideal
   $\mathfrak{I}'=\mathfrak{I}_1+\cdots+\mathfrak{I}_{r-1}$ in
   $\mathfrak{S}'$.  Now, combining with Proposition
   \ref{kringgenfbmind} we have
   \[K^0(M_r)\simeq \mathfrak{S}'/\mathfrak{I}'\bigotimes
     \mathbb{Z}[{\bf y}_{r}]^{W_{L_r}}\] modulo the ideal generated by the
   elements $g({\bf y}_r)-\bar{\Psi_r^*(g({\bf y}_r))}$ for
   $g({\bf y}_r)\in \mathbb{Z}[{\bf y}_{r}]^{W_r}$. Here
   $\bar{\Psi_r^*(g({\bf y}_r))}$ denotes the projection of
   $\Psi_r^*(g({\bf y}_r))\in \mathfrak{S}'$ in
   $\mathfrak{S}'/\mathfrak{I}'\simeq K^0(M_{r-1})\simeq
   \mathscr{K}^0(M_{r-1})$. This implies \eqref{presgenfbm}.  Hence the
   theorem.
 \end{proof}

 \subsection{Flag Bott manifold of type $A$}
 In this section we shall derive the results in \cite{pu} already
 obtained for the ordinary $K$-ring of the flag Bott manifolds of type
 $A$. The results in this section are $K$-theoretic analogue of the
 results on the equivariant cohomology in \cite[Corollary 4.4]{kkls}

 In this case $G_j=SL(n_j+1)$ and $P_j=B_j$ is the subgroup of upper
 triangular matrices with determinant $1$ for every $1\leq j\leq r$
 and $T_j$ is isomorphic to the subgroup of complex diagonal matrices
 $D_j\simeq (\mathbb{C}^*)^{n_j+1}$ with determinant $1$. Thus
 $K_j=SU(n_j+1)$ and $Z_j=(T_j)_{comp}$ is the maximal torus of
 $SU(n_j+1)$ for every $1\leq j\leq r$. Thus
 $K_j/Z_j=Flag(\mathbb{C}^{n_j+1})$ for each $1\leq j\leq r$. Thus we
 get the $r$-stage flag Bott tower constructed in \cite[Definition
 2.1]{klss}. In this case $W_j=S_{n_j+1}$ the symmetric group of
 $n_j+1$ elements for $1\leq j\leq r$ and thus $R(T_j)^{W_{j}}$ is
 generated by the elementary symmetric polynomials $e_{k}(\bf{y}_j)$
 in the variables
 $\{y_{j,1},\ldots, y_{j,n_j+1}\}$.  Moreover, since
 $L_j=T_j$ for $1\leq j\leq
 r$, the morphism of algebraic groups $\displaystyle\Psi_j:
 \prod_{l=1}^{j-1} L_l\lra
 T_{j}$ is equivalent to giving a collection of matrices
 $A_l^{(j)}=(A_l^{(j)}(i,s))\in M_{n_j+1,n_l+1}(\bz)$ for $1\leq i\leq
 n_j+1$, $1\leq s\leq n_l+1$ and $1\leq l\leq j-1$.

Let ${\bf{y}}_{j}$ denote the collection of variables
$\{y^{\pm1}_{j,1},\,y^{\pm 1}_{j,2},\ldots,y^{\pm1}_{j,n_j+1}\}$ for every $1\leq j\leq
r$.

Let $\mathcal{R}:=\bz[{\bf{y}}_{j}\,|\, 1\leq j\leq r]$ denote the
Laurent polynomial ring in the variables
$y^{\pm1}_{j,1},\,y^{\pm1}_{j,2},\ldots,y^{\pm1}_{j,n_j+1}$ for
$1\leq j\leq r$.

Let $e_{k}(\bf{y}_{j})\in \mathcal{R}$ denote the $k$th elementary
symmetric function in $y_{j,1},\,y_{j,2},\ldots,y_{j,n_j+1}$ for
$1\leq k\leq n_j+1$.

Let $\mathfrak{A}_l^{(j)}=(A_l^{(j)}(i,s))\in M_{n_j+1,n_l+1}(\bz)$ for
$1\leq i\leq n_j+1$, $1\leq s\leq n_l+1$.

Let $\mathcal{I}_1$ denote the ideal in $\mathcal{R}$ generated by the
polynomials
$$e_k({\bf{y}}_{1})- \begin{pmatrix} n_1+1\\k\end{pmatrix}, \text{ for
  every $1\leq k\leq n_1+1$}$$ and let $\mathcal{I}_j$ denote the
ideal in $\mathcal{R}$ generated by the polynomials
\[e_k({\bf{y}}_{j})-e_k\big(\prod_{l=1}^{j-1}(\prod_{s=1}^{n_l+1}y_{l,s}^{A^{(j)}_l(i,s)})\,|
  1\leq i\leq n_j+1\big)\] for every $2\leq j\leq r$ and
$1\leq k\leq n_j+1$. 

Let $\mathcal{I}:=\mathcal{I}_1+\cdots+\mathcal{I}_r$ which is the
smallest ideal in $\mathcal{R}$ containing
$\mathcal{I}_1,\ldots, \mathcal{I}_r$.

Thus when $G_j=SL(n_j+1)$ and $P_j=B_j$ the presentation in 
\eqref{presgenfbm} reduces to the following (see \cite{pu}).

\bth\label{prestypeAfbm}  
  Suppose that $M_\bullet =\{M_j\,|\,0\leq j\leq r\}$ is an $r$- stage
  flag Bott manifold of Lie type $A$, determined by a set of integer
  matrices
  $\mathfrak{A}:=(A_l^{(j)})_{1\leq l<j\leq r}\in \prod_{1\leq l<j\leq
    r}M_{n_j+1,n_l+1}(\bz)$. The map from $\mathcal{R}$ to
  $K^0(M_r)\simeq \mathscr{K}^0(M_{r})$
  which sends $y_{r,k}$ to $[\mathcal{L}_{r,k}]$ for
  $1\leq k\leq n_r+1$ and $y_{l,k}$ to
  \[[p_r^*\,\circ\, \cdots\,\circ\, p^*_{l+1}(\mathcal{L}_{l,k})]\] for
  $1\leq k\leq n_l+1$ and $1\leq l\leq r-1$,  induces
  an isomorphism $$\mathcal{R}/\mathcal{I}\simeq K^0(M_r)\simeq \mathscr{K}^0(M_r).$$
\eeth

\brem\label{linebundles} Recall that we have complex algebraic line
bundles $\mathcal{L}_{l,k}$ for $1\leq k\leq n_l+1$ on $M_l$ for $0\leq
l\leq r$ (see \cite{klss}  and \cite{pu}).
\erem

\subsection{Full flag Bott tower of general Lie type}

When $P_j=B_j$ and $L_j=T_j$ for every $1\leq j\leq r$, then $M_{j}$
is a bundle over $M_{j-1}$ with fiber the full flag $G_{j}/B_{j}$ for
$1\leq j\leq r$. In this case we shall denote
by \[M_{\bullet}^{full}:=\{M_j^{full}:=M_j~:~1\leq j\leq r\}\] and
call it the {\it full flag Bott tower of general Lie type}.

Let
\[\mathfrak{S}:=\bigotimes_{j=1}^r\mathbb{Z}[{\bf y}_j]\simeq
\bigotimes_{j=1}^r R(T_j).\] Let $\mathfrak{I}_j$
denote the ideal in $\mathfrak{S}$ generated by the elements
\[g({\bf y}_j)-g(\Psi_j^*({\bf y}_j))\] where
$g({\bf y}_j)\in R(T_j)^{W_j}=\mathbb{Z}[{\bf y}_j]^{W_j}$ for every $1\leq j\leq r$. 
(Note that
\[\Psi_j^*:R(T_j)=\mathbb{Z}[{\bf y}_j]\lra \bigotimes_{k=1}^{j-1} R(T_k)\simeq \bigotimes_{k=1}^{j-1} \mathbb{Z}[{\bf
	y}_{k}]\] is a $\mathbb{Z}$-algebra homomorphism.) Let
$\mathfrak{I}=\mathfrak{I}_1+\cdots+\mathfrak{I}_r$.

The following corollary of Theorem \ref{kringgenfbm} describes the
topological $K$-ring and the Grothendieck ring of $M_r^{full}$.

\begin{coro}\label{fullflagBott}
	 We have the
	following isomorphism of rings
	\be\label{equivpresgenfbm}\mathscr{K}^0(M_r^{full})\simeq
	K^0(M_r^{full})\simeq
	\mathfrak{S}/\mathfrak{I}.\ee
\end{coro}

\brem Since by \cite[Proposition 4.6, Corollary 4.7]{fls} a flag Bott
Samelson variety $Z$ of general Lie type is diffeomorphic to a full
flag Bott tower $X$ of general Lie type, Corollary \ref{fullflagBott}
gives the description of the $K$-ring of $Z$
generalizing the corresponding results for type $A$ in \cite[Section
4.2]{pu}.

\erem

\section{Equivariant $K$-ring of flag Bott manifolds of general Lie type}

\subsection{The torus action on flag Bott manifolds}\label{torusaction}

Let $\mathbb{T}$ denote the torus $\displaystyle\prod_{l=1}^r
T_l$. Then we have a well defined action of $\mathbb{T}$ on $M_j$
given by
\[(t_1,\ldots, t_r)\cdot [g_1,\ldots, g_j]:=[t_1\cdot g_1,\ldots,
  t_j\cdot g_j],\] where $(t_1,\ldots, t_r)\in \mathbb{T}$ and
$[g_1,\ldots, g_j]\in M_j$. The well-definedness follows by the fact
that the image of $\Psi^{(k)}_j$ lies in the commutative groups $T_l$ for
every $1\leq j<k\leq r$. (See \cite[Definition 4.1]{kkls} and
\cite{klss} where $\mathbb{D}$ denoted the above torus and
$\mathbb{T}$ denoted the maximal compact subgroup of the
$\mathbb{D}$.) Also note that the factors $T_l$ for $j<l\leq r$ act
trivially on $M_j$.

In this section we shall describe the
$\mathbb{T}_{comp}:=\prod_{j=1}^r (T_{j})_{comp}$-equivariant $K$-ring
(respectively $\mathbb{T}$-equivariant Grothendeick ring) and give its
presentation as an $R(\mathbb{T}_{comp})$-algebra (respectively
$R(\mathbb{T})$-algebra).

Let $H$ and $S$ be complex linear algebraic groups. Let ${E}\lra {Y}$
be an $S$-equivariant principal $H$-bundle where $H$ acts freely on
$E$ from right so that $Y=E/H$ and the group $S$ acts on $E$ and $Y$
on the left such that the bundle projection $E\lra E/H=Y$ is
$S$-equivariant. Let $\Phi:H\lra T$ be a morphism of algebraic groups
where $T:=(\mathbb{C}^*)^n$.

We then have the following $S$-equivariant analogue of Theorem
\ref{algcellularkunneth}.

\bth\label{equivalgcellularkunneth} For a $T$-cellular variety $X$ we
have the following isomorphism
\be\label{equivalgkunneth}\mathscr{K}_{S}^0({E}(X))\simeq
\mathscr{K}_{S}^0(Y)\otimes_{R(T)} \mathscr{K}^0_{T}(X)\ee (see
\cite[Corollary 1.8]{u}) where $\mathscr{K}_{S}^0$ denotes the
$S$-equivariant Grothendieck ring and $\mathscr{K}_{T}^0$ denotes the
$T$-equivariant Grothendieck ring. Here $\mathscr{K}_{S}^0(Y)$ gets
the structure of an $R(T)$-module by sending $[\mathbf{V}]\in R(T)$ to
\[[E\times_{H} \mathbf{V}]\in \mathscr{K}_{S}^0(Y),\] where a $T$-representation
$\mathbf{V}$ is thought of as a $H$-representation via the morphism
$\Phi$. Note that since $E$ has an $H$ action on the right and an
$S$-action on the left such that $E\lra E/H=Y$ is an $S$-equivariant
map, the vector bundle $[E\times_{H} \mathbf{V}]$ over $Y$ is
$S$-equivariant. \eeth
\begin{proof}
We consider the following map analogous to
$\phi'$ defined in \eqref{map1}.

\be\label{map1Y} \phi'_{S}:\mathscr{K}_{S}^0(Y)\otimes
\mathscr{K}^0_{T}(X)\lra \mathscr{K}_{S}^0({E}(X))\ee which takes
$[\mathcal{V}]\otimes [\mathcal{W}]$ to
$[\pi^*(\mathcal{V})]\otimes [E\times_{H} \mathcal{W}]$. Further, we
note that $\phi'$ induces the following map analogous to \eqref{map2}
\be\label{map2Y}\phi_{S}:\mathscr{K}_{S}^0(Y)\otimes_{R(T)}\mathscr{K}^0_{T}(X)\lra
\mathscr{K}^0_{S}({E}(X)).\ee The proof now follows verbatim as that
of Theorem \ref{algcellularkunneth}.
  \end{proof}

  By considering the maximal compact subgroup $S_{comp}$ of $S$ we
  have the following $S_{comp}$-equivariant analogue of Theorem
  \ref{cellularkunnethtop}.

  \bth\label{equivcellularkunnethtop} We have the following
  isomorphism: \be\label{equivtopkunneth}{K}_{S_{comp}}^0(E(X))\simeq
  {K}_{S_{comp}}^0(Y)\otimes_{R(T_{comp})}{K}^0_{T_{comp}}(X).\ee Here
  ${K}^0(Y)$ is a $R(T_{comp})=R(T)$-module via the map which takes
  any $[\mathbf{V}]\in R(T_{comp})$ for a $T$-representation
  $\mathbf{V}$ to $[E\times_{H} \mathbf{V}]\in {K}_{S_{comp}}^0(Y)$,
  where $H$ acts on $\mathbf{V}$ via $\Phi$.  \eeth

  By immediate application of Theorem \ref{equivalgcellularkunneth} to
  the $\mathbb{T}$-equivariant principal
  $P_1\times\cdots\times P_{r-1}$-bundle
  \[G_1\times \cdots\times G_{r-1}\lra (G_1\times \cdots\times
    G_{r-1})/(P_1\times\cdots\times P_{r-1})=M_{r-1}\] and the
  $T_r$-cellular variety $X=G_r/P_r$ we get the following
  $\mathbb{T}$-equivariant analogue of Proposition
  \ref{kringgenfbmind}.

\bpropo\label{equivkringgenfbmind}
\[ \mathscr{K}_{\mathbb{T}}^0(M_r)\simeq
  \mathscr{K}_{\mathbb{T}}^0(M_{r-1})
  \otimes_{R(T_r)^{W_r}} R(T_r)^{W_{L_{r}}}.\] Here
$\mathscr{K}_{\mathbb{T}}^0(M_{r-1})$ is a
$R(T_r)^{W_r}$-algebra via the sequence of maps
\[ \displaystyle R(T_r)^{W_r}\hra R(T_r)\stackrel{\Psi_r^*}{\lra}
  \bigotimes_{j=1}^{r-1} R(L_j)=R(\prod_{j=1}^{r-1} L_j)\] where
$\displaystyle[\mathbf{V}]\in R(T_r)^{W_r}$ maps to the class of
the $\mathbb{T}$-equivariant associated vector bundle
$\displaystyle\prod_{j=1}^{r-1}G_j\times_{\prod_{j=1}^{r-1} P_j}
\mathbf{V}\in \mathscr{K}_{\mathbb{T}}^0(M_{r-1})$ by considering
$\mathbf{V}$ as a
$\displaystyle\prod_{j=1}^{r-1} P_j$-representation via the maps
\[ \prod_{j=1}^{r-1} P_j\stackrel{\Phi_{r-1}}{\lra} \prod_{j=1}^{r-1}
  L_j\stackrel{\Psi_{r}}{\lra} T_r .\]  \epropo

Further, note that the factor $T_r$ acts trivially on $M_{r-1}$. Thus
\be\label{splitiso}\mathscr{K}^0_{\mathbb{T}}(M_{r-1})\simeq
\mathscr{K}^0_{\mathbb{T}'}(M_{r-1})\otimes R(T_r)\ee (see
\cite[(5.2.4)]{cg}) where
$\displaystyle\mathbb{T}':=\prod_{j=1}^{r-1} T_j$. The above
isomorphism is an isomorphism of $R(T_r)$-algebras and hence as
$R(T_r)^{W_r}$-algebras.

Consider a $\mathbb{T}$-equivariant line bundle $\mathcal{L}$ on
$M_{r-1}$. Since $T_r$-action on $M_{r-1}$ is trivial and since
$\mathcal{L}$ is locally trivial on $M_{r-1}$ and corresponds locally
to a character of $T_r$ and since $M_{r-1}$ is an irreducible variety, it
follows that as a $T_r$-equivariant line bundle $\mathcal{L}$ is a
trivial line bundle associated to a character $\chi$. Let $e^{\chi}$
denote the trivial line bundle on $M_{r-1}$ associated to the
character $\chi:\mathbb{T}\lra T_r\stackrel{\chi}{\lra}
\mathbb{C}^*$. Hence $\mathcal{L}\otimes e^{-\chi}$ is a line bundle
on $M_{r-1}$ with a trivial action of $T_r$. Therefore
$\mathcal{L}':=\mathcal{L}\otimes e^{-\chi}$ is a
$\mathbb{T}'$-equivariant line bundle on $M_{r-1}$. Thus it follows
that under \eqref{splitiso}
$[\mathcal{L}]\in \mathscr{K}^0_{\mathbb{T}}(M_{r-1})$ maps to
$[e^{\chi}]\otimes [\mathcal{L}'] \in R(T_r)\otimes
\mathscr{K}^0_{\mathbb{T}'}(M_{r-1})$.

Recall that
$R(T_j)=\mathbb{Z}[X^*(T_j)]\simeq\mathbb{Z}[{\bf y}_j]$ where
${\bf y}_j:=\{y_{j,1}^{\pm 1},\ldots, y^{\pm 1}_{j,n_j+1}\}$ for
$1\leq j\leq r$ (see \cite{st}).

Let ${\bf u}_j:=\{u_{j,1}^{\pm 1},\ldots, u^{\pm 1}_{j,n_j+1}\}$ for
$1\leq j\leq r$.

Then $R(L_j)=\mathbb{Z}[{\bf y}_j]^{W_{L_j}}=R(T_j)^{{W_{L_j}}}$ for
$1\leq j\leq r$ (see \cite{kk}).

In particular, under the above isomorphism \eqref{splitiso} an element
$y_{r,k}\in R(T_{r})$ maps to $u_{{r},k}\otimes \Psi_r^*(y_{r,k})$ for
$1\leq k\leq n_{r}+1$. Here we consider \[\Psi_j^*:R(T_j)\lra
  \prod_{k=1}^{j-1}\big(R(L_k)=R(T_k)^{W_{L_k}}\big)\subseteq
  \prod_{k=1}^{j-1} R(T_k).\]

We further let ${\bf u}_j\otimes \Psi_j^*({\bf y}_j)$ denote the
collection \be\{u_{j,1}\otimes \Psi_j^*(y_{j,1}),\ldots, u_{j,
  n_j+1}\otimes \Psi_j^*(y_{j, n_j+1}), u_{j,1}^{-1}\otimes
\Psi_j^*(y_{j,1})^{-1},\ldots, u_{j, n_j+1}^{-1}\otimes \Psi_j^*(y_{j,
  n_j+1})^{-1}\} \ee of elements in
$\displaystyle R(T_j)\otimes \prod_{k=1}^{j-1} R(T_k)^{W_{L_k}}$ where
  $\Psi_j^*(y_{j,k}^{-1})=\Psi_j^*(y_{j,k}) ^{-1}$ for
  $1\leq k\leq n_j+1$.

Let
 \[\mathfrak{S}:=\mathbb{Z}[{\bf u}_j, 1\leq j\leq r]\otimes \bigotimes_{j=1}^r\mathbb{Z}[{\bf y}_j]^{W_{L_j}}\simeq
   \bigotimes_{j=1}^r R(T_j)\otimes \bigotimes_{j=1}^r R(L_j).\] Let $\mathfrak{I}_j$
 denote the ideal in $\mathfrak{S}$ generated by the elements
 \[g({\bf y}_j)-g({\bf u}_j\otimes \Psi_j^*({\bf y}_j))\] where
 $g({\bf y}_j)\in R(T_j)^{W_j}=\mathbb{Z}[{\bf y}_j]^{W_j}\subseteq
 \mathbb{Z}[{\bf y}_j]^{W_{L_j}}$ for every $1\leq j\leq r$. 
(Note that
 \[\Psi_j^*:R(T_j)=\mathbb{Z}[{\bf y}_j]\lra \bigotimes_{k=1}^{j-1}\mathbb{Z}[{\bf
     y}_{k}]^{W_{L_k}}\subseteq \bigotimes_{k=1}^{j-1}\mathbb{Z}[{\bf
     y}_{k}]\] is a $\mathbb{Z}$-algebra homomorphism.) Also by
 \cite{st} $R(G_j)=R(T_j)^{W_j}$.  Let
 $\mathfrak{I}=\mathfrak{I}_1+\cdots+\mathfrak{I}_r$.

 The following theorem gives a presentation for the
 $\mathbb{T}_{comp}$-equivariant toplogical $K$-ring of $M_r$ and the
 $\mathbb{T}$-equivariant Grothendieck ring of $M_r$. See
 \cite[Corollary 4.3]{kkls} for the analogous statement for
 equivariant cohomology ring. \bth\label{equivkringgenfbm} We have the
 following isomorphism of rings
 \be\label{equivpresgenfbm}\mathscr{K}_{\mathbb{T}}^0(M_r)\simeq
 K^0_{\mathbb{T}_{comp}}(M_r)\simeq
 \mathfrak{S}/\mathfrak{I}.\ee\eeth \begin{proof} The first
   isomorphism follows by \cite[Proposition 5.5.6]{cg} since $M_r$ is
   constructed as iterated cellular fibrations. We now prove the
   second isomorphism in \eqref{equivpresgenfbm}. When $r=1$ then
   $\mathfrak{S}=\mathbb{Z}[{\bf u}_1]\otimes\mathbb{Z}[{\bf
     y}_1]^{W_{L_1}}$ and $\mathfrak{I}=\mathfrak{I}_1$ is generated
   by elements
   $g({\bf y}_1)-g({\bf u}_1\otimes\Psi_1^*({\bf y}_1))$ for
   $g({\bf y}_1)\in R(T_1)^{W_1}$ and $\Psi_1: pt=T_0\lra T_1$. Thus
   $\Psi_1^*: R(T_1)\lra \mathbb{Z}$ is the augmentation map which
   sends $[V]$ to $\mbox{dim}(V)$. But it is known that we have an
   isomorphism
   \[K^0_{T_1}(G_1/P_1)\simeq
   R(T_1)\otimes_{R(T_1)^{W_1}}R(T_1)^{W_{L_1}}\] as an
   $R(T_1)$-algebra. Thus the statement is true for $r=1$. We assume
   by induction that the statement is true for $r-1$. Thus
  $K^0_{\mathbb{T}_{comp}'}(M_{r-1})$ is isomorphic to
   \[\mathfrak{S}'= \mathbb{Z}[{\bf u}_j, 1\leq j\leq {r-1}]\otimes
     \bigotimes_{j=1}^{r-1}\mathbb{Z}[{\bf y}_j]^{W_{L_j}}\] modulo
   the ideal $\mathfrak{I}'=\mathfrak{I}_1+\cdots+\mathfrak{I}_{r-1}$
   in $\mathfrak{S}'$. This further implies that
   $K^0_{\mathbb{T}_{comp}}(M_{r-1})\simeq R(T_r)\otimes
   K^0_{\mathbb{T}_{comp}'}(M_{r-1})$ is isomorphic to
   \[\mathfrak{S}''= \mathbb{Z}[{\bf u}_j, 1\leq j\leq {r}]\otimes
     \bigotimes_{j=1}^{r-1}\mathbb{Z}[{\bf y}_j]^{W_{L_j}}\]
   modulo the ideal
   $\mathfrak{I}''=\mathfrak{I}_1+\cdots+\mathfrak{I}_{r-1}$ in
   $\mathfrak{S}''$.

   Now, combining with Proposition
   \ref{kringgenfbmind} we have $K_{\mathbb{T}_{comp}}^0(M_r)$ is isomorphic
   to
   \[\mathfrak{S}''/\mathfrak{I}''\bigotimes \mathbb{Z}[{\bf
       y}_{r}]^{W_{L_r}}\] modulo the ideal generated by the elements
   \[g({\bf y}_r)-g({\bf u}_r\otimes \bar{\Psi_r^*({\bf y}_r)})\] for
   $g({\bf y}_r)\in \mathbb{Z}[{\bf y}_{r}]^{W_r}$. Here
   ${\bf u}_r\otimes \bar{\Psi_r^*({\bf y}_r)}$ denotes the collection
   of elements \be\{u_{r,1}\otimes \bar{\Psi_r^*(y_{r,1})},\ldots,
   u_{r,n_r+1}\otimes \bar{\Psi_r^*(y_{r,n_r+1})}, u_{r,1}^{-1}\otimes
   \bar{{\Psi_r^*(y_{r,1})}}^{-1},\ldots u_{r,n_r+1}^{-1}\otimes
     \bar{\Psi_r^*(y_{r,n_r+1})}^{-1} \} \ee where
       $u_{r,k}\otimes \bar{\Psi_r^*(y_{r,n_k})}$ (respectively
       $u_{r,k}^{-1}\otimes \bar{\Psi_r^*(y_{r,n_k})}^{-1}$ ) denotes
       the projection of $u_{r,k}\otimes \Psi_r^*(y_{r,k})$
       (respectively $u_{r,k}^{-1}\otimes \Psi_r^*(y_{r,k})^{-1}$)
       in $\mathfrak{S}''$ in
       $\mathfrak{S}''/\mathfrak{I}''\simeq
       K^0_{\mathbb{T}_{comp}}(M_{r-1})\simeq
       \mathscr{K}^0_{\mathbb{T}}(M_{r-1})$. This implies \eqref{equivpresgenfbm}.
       Hence the theorem.
 \end{proof}

 \subsection{Equivariant $K$-ring of Flag Bott manifold of type $A$}
 
 Let ${\bf{y}}_{j}$ denote the collection of variables
 $\{y^{\pm1}_{j,1},\,y^{\pm 1}_{j,2},\ldots,y^{\pm1}_{j,n_j+1}\}$ for
 every $1\leq j\leq r$.

Let ${\bf u}_j:=\{u_{j,1}^{\pm 1},\ldots, u^{\pm 1}_{j,n_j+1}\}$ for
$1\leq j\leq r$.

Let
$\mathcal{S}:=\bz[{\bf{u}}_{j}\,|\, 1\leq j\leq r]\otimes
\bz[{\bf{y}}_{j}\,|\, 1\leq j\leq r]$. 

Let $e_{k}(\bf{y}_{j})\in \mathcal{S}$ denote the $k$th elementary
symmetric function in $y_{j,1},\,y_{j,2},\ldots,y_{j,n_j+1}$.

Let $\mathfrak{A}_l^{(j)}=(A_l^{(j)}(p,q))\in M_{n_j+1,n_l+1}(\bz)$ for
$1\leq p\leq n_j+1$, $1\leq q\leq n_l+1$.

Let $\mathcal{I}_1$ denote the ideal in $\mathcal{S}$ generated by the
polynomials
$$e_k({\bf{y}}_{1})-e_k({\bf u}_{1}), \text{ for
  every $1\leq k\leq n_1+1$}$$ and let $\mathcal{I}_j$ denote the
ideal in $\mathcal{S}$ generated by the polynomials
\[e_k({\bf{y}}_{j})-e_k(\big\{u_{j,i}\otimes
\prod_{l=1}^{j-1}\prod_{t=1}^{n_l+1}y_{l,t}^{A^{(j)}_l(i,t)},
u_{j,i}^{-1}\otimes \big(\prod_{l=1}^{j-1}\prod_{t=1}^{n_l+1}y_{l,t}^{A^{(j)}_l(i,t)}\big) ^{-1} | 1\leq
i\leq n_j+1\big\})\] for every $2\leq j\leq r$ and
$1\leq k\leq n_j+1$.



Let $\mathcal{I}:=\mathcal{I}_1+\cdots+\mathcal{I}_r$ which is the
smallest ideal in $\mathcal{S}$ containing
$\mathcal{I}_1,\ldots, \mathcal{I}_r$.

Thus when $G_j=SL(n_j+1)$ and $P_j=B_j$ the presentation in 
\eqref{equivpresgenfbm} reduces to the following presentation of
$K^0_{\mathbb{T}_{comp}}(M_{r})$.

\bth\label{equivprestypeAfbm} Suppose that
$M_\bullet =\{M_j\,|\,0\leq j\leq r\}$ is an $r$- stage flag Bott
manifold of Lie type $A$, determined by a set of integer matrices
$\mathfrak{A}:=(A_l^{(j)})_{1\leq l<j\leq r}\in \prod_{1\leq l<j\leq
  r}M_{n_j+1,n_l+1}(\bz)$. The map from $\mathcal{S}$ to
$K_{\mathbb{T}_{comp}}^0(M_r)$ which sends $y_{r,k}$ to
$[\mathcal{L}_{r,k}]_{\mathbb{T}_{comp}}$ for $1\leq k\leq n_r+1$ and
$y_{l,k}$ to
  \[[p_r^*\,\circ\, \cdots\,\circ\,
    p^*_{l+1}(\mathcal{L}_{l,k})]_{\mathbb{T}_{comp}}\] for
  $1\leq k\leq n_l+1$ and $1\leq l\leq r-1$, induces an isomorphism
  $$\mathcal{S}/\mathcal{I}\simeq K_{\mathbb{T}_{comp}}^0(M_r).$$ \eeth

\subsubsection{Bott manifold of type A}
Suppose that an flag Bott tower $M^{Bott}_\bullet:=\{M^{Bott}_j: j=1,2,\ldots,r\}$ is an r- stage
Bott tower of Lie type A, determined by a set of integer matrices
$\mathfrak{A}:=(A^{(j)}_l)_{1\leq l<j\leq r}\in M_{2\times 2}(\bz)$ as each $n_j=1$ for $ j=1,2,\ldots,r$.

We let $c_{j,1}$ and $c_{j,2}$ to be \[u_{j,1} \otimes
	\prod_{l=1}^{j-1}(y_{l,1}^{A^{(j)}_l(1,1)}\cdot y_{l,2}^{A^{(j)}_l(1,2)}) \text { and } u_{j,2} \otimes
	\prod_{l=1}^{j-1}(y_{l,1}^{A^{(j)}_l(2,1)}\cdot y_{l,2}^{A^{(j)}_l(2,2)})\] respectively for all $2\leq j\leq r$. Also, we let $c_{1,1}=u_{1,1}$ and $c_{1,2}=u_{1,2}$ when $j=1$.

Let $\mathcal{S}\simeq \bz[y_{j,1}^{\pm1},\,y_{j,2}^{\pm1},\,u_{j,1}^{\pm1},\,u_{j,2}^{\pm1}|1\leq j\leq r]$.

 Let $\mathcal{I}_j$ for $1\leq j\leq r$ denote the
ideal in $\mathcal{S}$ generated by the polynomials
\[(y_{j,1}+y_{j,2})-(c_{j,1}+c_{j,2}) \text{ and } (y_{j,1}\cdot y_{j,2})-(c_{j,1}\cdot c_{j,2})\]

Let $\mathcal{I}:=\mathcal{I}_1+\cdots+\mathcal{I}_r$ which is the
smallest ideal in $\mathcal{S}$ containing
$\mathcal{I}_1,\ldots, \mathcal{I}_r$.  The following corollary of
Theorem \ref{equivprestypeAfbm} gives the presentation of
$T_{comp}$-equivariant $K$-ring of a Bott manifold.

\begin{coro}
	 We have the
	following isomorphism of rings
	\be\label{equivpresBott}\mathscr{K}_{\mathbb{T}}^0(M_r^{Bott})\simeq
	K^0_{\mathbb{T}_{comp}}(M_r^{Bott})\simeq
	\mathfrak{S}/\mathfrak{I}.\ee
\end{coro}
\begin{rema}
  From \eqref{equivpresBott}, we can also describe the topological
  equivariant $K$ ring of a Bott Samelson variety as it degenerates to
  a Bott manifold \cite{gk} preserving the underlying differentiable
  structure. This has also been done earlier in \cite[Theorem 5.7,
  Theorem 6.3]{w2} by alternate methods.
\end{rema}

\section{Examples}  

The following example illustrates the computation of the
  $\mathbb{T}_{comp}$-equivariant $K$-ring of a flag Bott manifold of
  height $r$ with all the fibres full flag manifolds $Sp(n)/T_n$ where $T_n$
  denotes the maximal torus of $Sp(n)$ (see \cite[p. 193]{hus}).
 
  \beg Consider $Sp(n)\subseteq U(2n)$ as matrices with quaternion
  entries by identifying $\mathbb{H}^n$ with $\mathbb{C}^{2n}$. The
  maximal torus $T_{n}$ of $Sp(n)$ are diagonal matrices with entries in
  the circle subgroup $S^1\subseteq \mathbb{C}\subseteq \mathbb{H}$.

  Let
  $M_r=Sp(n_1)\times_{T_{n_1}}\times\cdots\times
  Sp(n_{r-1})\times_{T_{n_{r-1}}}Sp(n_r)/T_{n_r}$.

  A morphism $(S^1)^{n_{l}}\simeq T_{n_l}\lra T_{n_j} \simeq
  (S^1)^{n_j}$ is determined by an integer matrix $A_l^{(j)}\in
  M_{n_j\times n_l}(\mathbb{Z})$. Thus $\displaystyle \Psi_{j}:\prod_{l=1}^{j-1}
  T_{n_l}\lra T_{n_j}$ is given by a collection of matrices
  $A^{(j)}_l\in M_{n_j\times n_l}(\mathbb{Z})$ for  $1\leq l\leq j-1$.

  In this case, $P_j=B_j$, $L_j=T_{n_j}$ for $1\leq j\leq r$ and
  $W_{L_j}=\{1\}$ for $1\leq j\leq r$. Thus
  $R(T_{n_j})=\mathbb{Z}[{u_{j,1}}^{\pm 1},\ldots, {u_{j,n_j}}^{\pm 1}]$ and
  $R(Sp(n_j))=R(T_{n_j})^{W_j}=\mathbb{Z}[\lambda_1,\ldots,
  \lambda_{n_j}]$ where $\lambda_k$ is the $k$th elementary symmetric
  function in the $2n_j$ variables
  $u_{j,1}, u^{-1}_{j,1},\ldots, u_{j,n_j}, u^{-1}_{j,n_j}$ (see \cite[p. 195]{hus}).

Let ${\bf{y}}_{j}$ denote the collection of variables
 $\{y^{\pm1}_{j,1},\,y^{\pm 1}_{j,2},\ldots,y^{\pm1}_{j,n_j}\}$ for
 every $1\leq j\leq r$.

Let ${\bf u}_j:=\{u_{j,1}^{\pm 1},\ldots, u^{\pm 1}_{j,n_j}\}$ for
$1\leq j\leq r$.
  
By Theorem \ref{equivkringgenfbm} we
have
\[\mathscr{K}_{\mathbb{T}}^0(M_r)\simeq
  K^0_{\mathbb{T}_{comp}}(M_r)\simeq \mathfrak{S}/\mathfrak{I} \]
where

 \[\mathfrak{S}:=\mathbb{Z}[{\bf u}_j, 1\leq j\leq r]\otimes \bigotimes_{j=1}^r\mathbb{Z}[{\bf y}_j]\simeq
   \bigotimes_{j=1}^r R(T_j)\otimes \bigotimes_{j=1}^r R(T_j)\] and
 $\mathfrak{I}_j$ denote the ideal in $\mathfrak{S}$ generated by the
 elements
 \[g({\bf y}_j)-g({\bf u}_j\otimes \Psi_j^*({\bf y}_j))\] where
 $g({\bf y}_j)\in \mathbb{Z}[{\bf y}_j]^{W_j}\simeq
 \mathbb{Z}[\lambda_1,\ldots, \lambda_{n_j}]\subseteq \mathbb{Z}[{\bf
   y}_j]$ for every $1\leq j\leq r$. Here $\lambda_1,\ldots,
 \lambda_{n_j}$ are the elementary symmetric functions in ${\bf y}_j$.
 \eeg

The following is an example of a full flag Bott tower of height $3$ of
mixed types $\mbox{A}$, $\mbox{C}$ and $\mbox{B}$.

\begin{exam}
	Let $T_{n_1}$ be a maximal torus of $SU(n_1)$, $T_{n_2}$ be a maximal torus of $Sp(n_2)$ and $W_1$ and $W_2$ be the respective Weyl group of them.
	
	Let ${\bf{e}}_1,\ldots,{\bf{e}}_{2n_3+1}$ denote the canonical basis of $\br^{2n_3+1}$. We define $w_j: S^1\simeq \br/\bz\to Spin(2n_3+1)$ to be the homomorphism given by the relation $w_j(\theta)=cos \,2\pi\theta \, + {\bf{e}}_{2j-1}\,{\bf{e}}_{2j} \,sin \,2\pi\theta$ for each $1\leq j\leq n_3$.  
	
	\noindent
	Let $w: (S^1)^{n_3} \to  Spin(2n_3+1)$ be defined by \[w(\theta_1,\ldots,\theta_{n_3})= w_1(\theta_1)\cdots w_{n_3}(\theta_{n_3}).\]
	
	\noindent
	Define  $T_{n_3}:= w((S^1)^{n_3})$ which is a maximal torus of $Spin(2n_3+1)$ and $W_3$ to be the Weyl group of $Spin(2n_3+1)$, which is isomorphic to the Weyl group of $SO(2n_3+1)$ (see \cite[Proposition 8.2, p. 197]{hus}) consisting of the $2^{n_3}n_3!$ permutations of the indexes of $(\theta_1,\ldots ,\theta_{n_3})$ composed with substitutions $(\theta_1,\ldots,\theta_{n_3})\to(\pm \theta_1,\ldots,\pm \theta_{n_3})$.

   Let $M_3$ be a 3 stage flag Bott manifold defined by \[M_3:=SU(n_1)\times_{T_{n_1}}Sp(n_2)\times_{T_{n_2}}Spin(2n_3+1)/T_{n_3}\] with homomorphisms $\Psi^{(j)}_l:T_{n_l}\to T_{n_j}$ which are determined by a collection of matrices $A^{(j)}_l\in M_{n_j\times n_l}(\bz)$ for $1\leq l< j\leq 3$.
	
	Let ${\bf{y}}_j$ (resp. ${\bf{u}}_j$) denote the collection of variables $\{y_{j,1}^{\pm1},\ldots,y_{j,n_j}^{\pm1}\}$ (resp. $\{u_{j,1}^{\pm1},\ldots,u_{j,n_j}^{\pm1}\}$) for each $j=1,2,3$.
	
     Note that $R(T_{n_1})=\bz[{\bf{u}}_1]$, $R(T_{n_2})=\bz[{\bf{u}}_2]$ and
    \[R(T_{n_3})=\bz[{\bf{u}}_3, (u_{3,1}\cdots u_{3,n_3})^{1/2}].\]

    $Spin(2n_3+1)$ is simply connected and twofold covering space of
    $SO(2n_3+1)$. Let $T_{n_3}'$ denote the maximal torus of
    $SO(2n_3+1)$. This is known  $ R(T_{n_3}')=\bz[{\bf{u}}_3]$. In
    case of $R(T_{n_3})$, there is one more extra character
    $(u_{3,1}\cdots u_{3,n_3})^{1/2}$ (see, \cite[p. 197, Proposition 8.3]{hus}). One may easily observe $(u_{3,1}\cdots u_{3,n_3})^{-1/2}\in R(T_{n_3})$ as \[(u_{3,1}\cdots u_{3,n_3})^{-1/2}=(u_{3,1}\cdots u_{3,n_3})^{1/2}\cdot(u_{3,1}^{-1}\cdots u_{3,n_3} ^{-1}).\]

    \noindent
    Let \[\mathfrak{S}=\bz[{\bf{y}}_1,\,{\bf{y}}_2,\,{\bf{y}}_3,\,(y_{3,1}\cdots y_{3,n_3})^{1/2}]\otimes \bz[{\bf{u}}_1,\,{\bf{u}}_2,\,{\bf{u}}_3, (u_{3,1}\cdots u_{3,n_3})^{1/2}].\]
	
	$\mathfrak{I}_1$ denote the ideal in $\mathfrak{S}$ generated by the elements
	\[e_k({\bf{y}}_1)-e_k({\bf{u}}_1)\] for $1\leq k\leq n_1$. Note that $e_k({\bf{y}}_1)\in R(T_{n_1})^{W_1}$ is the $k$-th elementary symmetric functions in the variables $\{y_{1,1},\ldots,y_{1,n_1}\}$.
	
	$\mathfrak{I}_2$ denote the ideal in $\mathfrak{S}$ generated by the elements
	\[g({\bf{y}}_2)-g({\bf{u}}_2\otimes \Psi_2^*({\bf{y}}_2)),\] where $g({\bf{y}}_2) \in R(T_{n_2})^{W_2}$ denote the $k$-th elementary symmetric functions in the variables $\{y_{2,1},\ldots,y_{2,n_2},y_{2,1}^{-1},\ldots,y_{2,n_2}^{-1}\}$ for $1\leq k\leq n_2$.
	
	and, $\mathfrak{I}_3$ denote the ideal in $\mathfrak{S}$ generated by the elements	\[g({\bf{y}}_3)-g({\bf{u}}_3\otimes \Psi_3^*({\bf{y}}_3)),\] where (see \cite[p. 200]{hus})
   \[\begin{split}
	g({\bf{y}}_3)\in  R(T_{n_3})^{W_3} & \simeq R(Spin(2n_3+1))\\
		&\simeq \bz[\lambda^1(\rho_{2n_3+1}),\ldots, \lambda^{n_3-1}(\rho_{2n_3+1}),\Delta_{2n_3+1}].
	\end{split}\]
     such that $\lambda^i(\rho_{2n_3+1})$ is the coeffecient of $t^i$ in the polynomial
	 \[(1+t)\cdot(1+y_{3,1}t)\cdot(1+y_{3,1}^{-1}t)\cdots(1+y_{3,n_3}t)\cdot(1+y_{3,n_3}^{-1}t).\]
	 and, \[\Delta_{2n_3+1}=\prod_{i=1}^{n_3} (y_{3,i}^{1/2}+y_{3,i}^{-1/2})=\sum_{\epsilon(l)=\pm 1} y_{3,1}^{\epsilon(1)/2}\cdots y_{3,n_j}^{\epsilon(n_j)/2}\]
	 Let $\mathfrak{I}$ be the smallest ideal in $\mathfrak{S}$ containing $\mathfrak{I}_1$, $\mathfrak{I}_2$ and $\mathfrak{I}_3$.
	 
	Hence, from Theorem \ref{equivkringgenfbm},
	 \[\mathscr{K}^0_\bt(M_3)\cong K_{\bt_{comp}}^0(M_3)\cong \mathfrak{S}/\mathfrak{I}.\]
\end{exam}

\vspace{0.5cm}
{\bf Acknowledgement:} The first named author
thanks Indian Institute of Technology, Madras for PhD fellowship.

\vspace{0.5cm}
\noindent{\bf Statements and Declarations} We hereby declare that this work
has no related financial or non-financial conflict of interests.

\end{document}